\documentclass[12pt]{amsart}

\usepackage{amssymb}
\usepackage{latexsym}
\usepackage{amsmath}

\begin{document}

\newtheorem{proposition}{Proposition}[section]
\newtheorem{definition}{Definition}[section]
\newtheorem{lemma}{Lemma}[section]
\newtheorem{theorem}{Theorem}[section]
\newtheorem{corollary}{Corollary}[section]
\newtheorem{problem}{Problem}[section]
\newtheorem{conjecture}{Conjecture}[section]
\newtheorem{remark}{Remark}[section]

\title{Subfactors and Hadamard Matrices}
\author{Wes Camp}
\address{W.C.: Department of Mathematics, Vanderbilt University, 1326 Stevenson Center, Nashville, TN 37240, USA}
\email{wes.camp@vanderbilt.edu}
\author{Remus Nicoara}
\address{R.N.: Department of Mathematics, Vanderbilt University, 1326 Stevenson Center, Nashville, TN 37240, USA}
\email{remus.nicoara@vanderbilt.edu}

\begin{abstract}
To any complex Hadamard matrix $H$ one associates a spin model commuting square, and therefore a hyperfinite subfactor. The standard invariant of this subfactor captures certain "group-like" symmetries of $H$. To gain some insight, we compute the first few relative commutants of such subfactors for Hadamard matrices of small dimensions. Also, we show that subfactors arising from Dita type matrices have intermediate subfactors, and thus their standard invariants have some extra structure besides the Jones projections. 
\end{abstract}

\maketitle

\section{Introduction} A complex Hadamard matrix is a matrix $H\in M_n(\mathbb C)$ having all entries of absolute value 1 and all rows mutually orthogonal. Equivalently, $\frac{1}{\sqrt{n}}H$ is a unitary matrix with all entries of the same absolute value. For example, the \emph{Fourier matrix} $F_n=(\omega^{ij})_{1\leq i,j\leq n}$, $\omega=e^{2\pi \mathfrak i/n}$, is a Hadamard matrix. 

In the recent years, complex Hadamard matrices have found applications in various topics of mathematics and physics, such as quantum information theory, error correcting codes, cyclic n-roots, spectral sets and Fuglede's conjecture. A general classification of real or complex Hadamard matrices is not available. A catalogue of most known complex Hadamard matrices can be found in \cite{TZ}. The complete classification is known for $n\leq 5$ (\cite{Haagerup}) and for self-adjoint matrices of order $6$ (\cite{BeN}).

The connection between Hadamard matrices and von Neumann algebras arose from an observation of Popa (\cite{Po2}): a unitary matrix $U$ is of the form $\frac{1}{\sqrt{n}}H$, $H$ Hadamard matrix, if and only if the algebra of $n\times n$ diagonal matrices $\mathcal D_n$ is orthogonal onto $U\mathcal D_n U^*$, with respect to the inner product given by the trace on $M_n(\mathbb C)$. Equivalently, the square of inclusions:
$$\mathfrak{C}(H)=\left(\begin{matrix}
\mathcal D_n & \subset{} & M_n(\mathbb{C}) \cr
\cup   &           & \cup \cr
\mathbb{C} & \subset{} & U\mathcal D_nU^*
\end{matrix},\tau\right)$$
is a \emph{commuting square}, in the sense of \cite{Po1},\cite{Po2}. Here $\tau$ denotes the trace on $M_n(\mathbb C)$, normalized such that $\tau(1)=1$. 

Such commuting squares are called \emph{spin models}, the name coming from statistical mechanical considerations (see \cite{JS}). By iterating Jones' basic construction, one can construct a hyperfinite, index $n$ subfactor from $H$ (see for instance \cite{JS}). The subfactor associated to $H$ can be used to capture some of the symmetries of $H$, and thus to classify $H$ to a certain extent (see \cite{BHJ},\cite{Jo2},\cite{BaN}). 

Let $N\subset M$ be an inclusion of $II_1$ factors of finite index, and let $N\subset{M}\overset{e_1}\subset{M_1}\overset{e_2}\subset M_2\subset{}...$ be the tower of factors constructed by iterating Jones' basic construction (see \cite{Jo1}), where $e_1,e_2,...$ denote the Jones projections. The \emph{standard invariant} $\mathcal G_{N,M}$ is then defined as the trace preserving isomorphism class of the following sequence of \emph{commuting squares} of inclusions of finite dimensional $*$-algebras:
$$
\begin{matrix}

\mathbb{C}=N'\cap N & \subset & N'\cap{M} & \subset & N'\cap{M_1} & \subset & N'\cap{M_2} & \subset &  ... & \cr
& & \cup    &          &  \cup            &          & \cup           &     &     & \cr
& & M'\cap{M} &\subset &M'\cap{M_1} & \subset & M'\cap{M_2} & \subset  & ... &  \cr 

\end{matrix}$$

The Jones projections $e_1,e_2,...,e_n$ are always contained in $N'\cap M_n$. If the index of the subfactor $N\subset M$ is at least $4$, they generate the Temperley-Lieb algebra of order $n$, denoted $TL_n$. In a lot of situations the relative commutant $N'\cap M_n$ has some interesting extra structure, besides $TL_n$. For instance, the five non-equivalent real Hadamard matrices of order $16$ yield different dimensions for the second relative commutant $N'\cap M_1$, and thus are classified by these dimensions (\cite{BHJ}).

In this paper we investigate the relation between Hadamard matrices and their subfactors. We look at Hadamard matrices of small dimensions or of special types. The paper is organized as follows: in the first section we recall, in our present framework, several results of \cite{Jo2},\cite{JS} regarding computations of standard invariants for spin models. 

In section 2 we study the subfactors associated to Hadamard matrices of Dita type. These are matrices that arise from a construction of \cite{Di}, which is a generalization of a construction of Haagerup (\cite{Haagerup}). Most known parametric families of Hadamard matrices are of Dita type. We show that the associated subfactors have intermediate subfactors. 

In the last section we present a list of computations of the second and third relative commutants $N'\cap M_1, N'\cap M_2$, for complex Hadamard matrices of small dimensions. We make several remarks and conjectures regarding the structure of the standard invariant. Most of the computations included were done using computers, with the help of the Mathematica and GAP softwares.

We would like to thank Teodor Banica, Kyle Beauchamp and Dietmar Bisch for fruitful discussions and correspondence. Wes Camp was supported in part by NSF under Grant No. DMS 0353640 (REU Grant), and Remus Nicoara was supported in part by NSF under Grant No. DMS 0500933.

\section{Subfactors associated to Hadamard matrices}
Let $H$ be a complex $n\times n$ Hadamard matrix and let $U=\frac{1}{\sqrt{n}}H$. $U$ is a unitary matrix, with all entries of the same absolute value. One associates to $U$ the square of inclusions:

$$\mathfrak{C}(H)=\left(\begin{matrix}
\mathcal D_n & \subset{} & M_n(\mathbb{C}) \cr
\cup   &           & \cup \cr
\mathbb{C} & \subset{} & U\mathcal D_nU^*
\end{matrix},\tau\right)$$
where $\mathcal D_n$ is the algebra of diagonal $n\times n$ matrices and $\tau$ is the trace on $M_n(\mathbb C)$, normalized such that $\tau(1)=1$.

 Since $H$ is a Hadamard matrix, $\mathfrak C(H)$ is a \emph{commuting square} in the sense of \cite{Po1},\cite{Po2}, i.e. $E_{\mathcal D_n}E_{U\mathcal D_nU^*}=E_{\mathbb C}$. The notation $E_A$ refers to the $\tau$-invariant conditional expectation from $M_n(\mathbb C)$ onto the $*$-subalgebra $A$.

Recall that two complex Hadamard matrices are said to be \emph{equivalent} if there exist unitary diagonal matrices $D_1,D_2$ and permutation matrices $P_1,P_2$ such that $H_2=P_1D_1H_1D_2P_2$. It is easy to see that $H_1,H_2$ are equivalent if and only if $\mathfrak C(H_1),\mathfrak C(H_2)$ are isomorphic as commuting squares, i.e. conjugate by a unitary from $M_n(\mathbb C)$.

We denote by $\mathfrak C^{t}(H)$ the commuting square obtained by flipping the upper left and lower right corners of $\mathfrak C(H)$:
$$\mathfrak{C}^{t}(H)=\left(\begin{matrix}
U\mathcal D_n U^* & \subset{} & M_n(\mathbb{C}) \cr
\cup   &           & \cup \cr
\mathbb{C} & \subset{} & \mathcal D_n
\end{matrix},\tau\right)$$

We have: $\mathfrak C^{t}(H)=\text{Ad}(U)\mathfrak C(H^*)$. Thus, $\mathfrak C^{t}(H)$ and $\mathfrak C(H)$ are isomorphic as commuting squares if and only if $H,H^*$ are equivalent as Hadamard matrices. 

We now recall the construction of a subfactor from a commuting square. By iterating Jones' basic construction (\cite{Jo1}), one obtains from $\mathfrak{C}^t(H)$ a tower of commuting squares of finite dimensional $*$-algebras:

\begin{equation}\label{tower}\begin{matrix}
U\mathcal D_n U^*& \subset{} & M_n(\mathbb{C}) & \overset{g_3}\subset & \mathcal X_1 & \overset{g_4}\subset & \mathcal X_2 & \overset{g_5}\subset &... \cr
\cup   &           & \cup & & \cup & &\cup & \cr
\mathbb{C} & \subset{} & \mathcal D_n & \overset{g_3}\subset & \mathcal Y_1 & \overset{g_4}\subset & \mathcal Y_2 & \overset{g_5}\subset &...
\end{matrix}\end{equation}
together with the extension of the trace, which we will still denote by $\tau$, and Jones projections $g_{i+2}\in \mathcal Y_i$, $i=1,2,...$. 

Let $M_H$ be the weak closure of $\cup_{i}\mathcal X_i$, with respect to the trace $\tau$, and let $N_H$ be the weak closure of $\cup_{i}\mathcal Y_i$. $N_H,M_H$ are hyperfinite $II_1$ factors, and the trace $\tau$ extends continuously to the trace of $M_H$, which we will still denote by $\tau$. It is well known that $N_H\subset M_H$ is a subfactor of index $n$, which we will call the subfactor associated to the Hadamard matrix $H$.

The standard invariant of $N_H\subset M_H$ can be expressed in terms of commutants of finite dimensional algebras, by using Ocneanu's compactness argument (5.7 in \cite{JS}). Consider the basic construction for the commuting square $\mathfrak C(H)$:
\begin{equation}\label{tower}\begin{matrix}
\mathcal D_n & \subset{} & M_n(\mathbb{C}) & \overset{e_3}\subset & \mathcal P_1 & \overset{e_4}\subset & \mathcal P_2 & \overset{e_5}\subset &... \cr
\cup   &           & \cup & & \cup & &\cup & \cr
\mathbb{C} & \subset{} & U\mathcal D_nU^* & \overset{e_3}\subset & \mathcal Q_1 & \overset{e_4}\subset & \mathcal Q_2 & \overset{e_5}\subset &...
\end{matrix}\end{equation}

Ocneanu's compactness theorem asserts that the first row of the standard invariant of $N_H\subset M_H$ is the row of inclusions:

$$\mathcal D_n'\cap U\mathcal D_nU^*\subset \mathcal D_n'\cap \mathcal Q_1 \subset \mathcal D_n'\cap \mathcal Q_2\subset \mathcal D_n'\cap \mathcal Q_3\subset ...$$

More precisely, if $$N_{H}\subset M_{H}\overset{e_3}\subset M_{H,1}\overset{e_4}\subset M_{H,2}\overset{e_5} \subset ...$$ is the Jones tower obtained from iterating the basic construction for the inclusion $N_H\subset M_H$, then:
$$D_n'\cap \mathcal Q_i=N_{H}'\cap M_{H,i},\text{ for all }i\geq 1.$$

Thus, the problem of computing the standard invariant of the subfactor associated to $H$ is the same as the computation of $\mathcal D_n'\cap \mathcal Q_i$. However, such computations seem very hard, and even for small $i$ and for matrices $H$ of small dimensions they seem to require computer use. Jones (\cite{Jo2}) provided a diagrammatic description of the relative commutants $\mathcal D_n'\cap \mathcal Q_i$ (see also \cite{JS}), which we express below in the framework of this paper.

Let $\mathcal P_0=M_n(\mathbb C)$ and let $(e_{i,j})_{1\leq i,j\leq n}$ be its canonical matrix units. Let $$e_2=\frac{1}{n}\sum_{i,j=1}^n e_{i,j}.$$

It is easy to check that $e_2$ is a projection. Moreover: $<\mathcal D_n, e_2>=M_n(\mathbb C)$ and $e_2xe_2=E_{\mathcal D_n}(x)e_2$ for all $x\in M_n(\mathbb C)$. Thus, $e_2$ is realizing the basic construction $$\mathbb C\subset \mathcal D_n\overset{e_2}\subset M_n(\mathbb C)$$

Let $e_{k,l}\otimes e_{i,j}$ denote the $n^2\times n^2$ matrix having only one non-zero entry, equal to $1$, at the intersection of row $(i-1)n+k$ and column $(j-1)n+l$. Thus, $e_{k,l}\otimes e_{i,j}$ are matrix units of $M_n(\mathbb C)\otimes M_n(\mathbb C)$. In what follows, we will assume that the embedding of $M_n(\mathbb C)$ into $M_n(\mathbb C)\otimes M_n(\mathbb C)$ is realized as $e_{k,l}\rightarrow e_{k,l}\otimes I_n$, where $e_{k,l}\otimes I_n=\sum_{i=1}^n e_{k,l}\otimes e_{i,i}$. 

\begin{lemma} Let $\mathcal P_1=M_n(\mathbb C)\otimes \mathcal D_n$, $\mathcal P_2=M_n(\mathbb C)\otimes M_n(\mathbb C)$, $e_3=\sum_{i=1}^n e_i\otimes e_i\in \mathcal P_1$ and $e_4=I_n\otimes e_2\in \mathcal P_2$. Then $$\mathcal D_n\subset M_n(\mathbb C)\overset{e_3}\subset \mathcal P_1$$ is a basic construction with Jones projection $e_3$ and $$M_n(\mathbb C)\subset \mathcal P_1\overset{e_4}\subset \mathcal P_2$$ is a basic construction with Jones projection $e_4$.
\end{lemma}
\begin{proof} To show that $\mathcal D_n\subset M_n(\mathbb C)\overset{e_3}\subset \mathcal P_1$ is a basic construction it is enough to check that $<M_n(\mathbb C), e_3>=\mathcal P_1$ and $e_3$ is implementing $E^{\mathcal P_1}_{M_n(\mathbb C)}$. First part is clear, since $e_{k,i}e_3e_{i,l}=e_{k,l}\otimes e_{i,i}$ are a basis for $\mathcal P_1=M_n(\mathbb C)\otimes \mathcal D_n$. To check that $e_3$ implements the conditional expectation, let $X=(x_{i,j})\in M_n(\mathbb C)$. We have:
\begin{equation}\begin{aligned}e_3(X\otimes I_n)e_3&=\sum_{i,j=1}^n(D_i\otimes D_i)(X\otimes I_n)(D_j\otimes D_j)\\
&=\sum_{i=1}^n D_iXD_i\otimes D_i\\
&=\sum_{i=1}^n (D_iXD_i\otimes I_n)e_3\\
&=E_{\mathcal D_n\otimes I_n}(X)e_3
\end{aligned}\end{equation}
 
Since $\mathbb C\subset \mathcal D_n\overset{e_2}\subset M_n(\mathbb C)$ is a basic construction, after tensoring to the left by $M_n(\mathbb C)$ it follows that $M_n(\mathbb C)\subset \mathcal P_1\overset{e_4}\subset \mathcal P_2$ is a basic construction, with $e_4=I_n\otimes e_2$.

 \end{proof}

\begin{proposition} The algebras $\mathcal P_1,\mathcal P_2,\mathcal P_3,...$ constructed in (\ref{tower}) are given by $$\mathcal P_{2k}=\otimes_{i=1}^{k+1}M_n(\mathbb C),\text{ } \mathcal P_{2k+1}=\mathcal P_{2k}\otimes \mathcal D_n$$ with the Jones projections $$e_{2k+2}=\otimes_{i=1}^k I_n\otimes e_2,\text{ }  e_{2k+3}=\otimes_{i=1}^k I_n\otimes e_3$$ 
\end{proposition}

\begin{proof} Follows from the previous lemma, by tensoring successively by $M_n(\mathbb C)$.
\end{proof}

\begin{proposition} Let $H$ be a complex $n\times n$ Hadamard matrix, let $U=\frac{1}{\sqrt{n}}H$, and $$D_U=\sum_{i,j=1}^n \bar u_{i,j} e_{j,j}\otimes e_{i,i}, \text{ }U_1=UD_U.$$
Then the algebras $\mathcal Q_1,\mathcal Q_2,\mathcal Q_3,...$ constructed in (\ref{tower}) are given by $$\mathcal Q_k=U_k\mathcal P_{k-1}U_k^*,\text{ }k\geq 1$$
where $U_k\in \mathcal P_k$ are the unitary elements: $$U_{2k+1}=\Pi_{i=0}^k(\otimes^i I_n\otimes U_1\otimes^{k-i} I_n),\text{ } U_{2k}=U_{2k-1}(\otimes^kI_n\otimes U),\text{ }k\geq 1.$$

\end{proposition}

\begin{proof} The unitary $U_1$ satisfies:
$$(\text{Ad}U_1)(\mathcal D_n)=(\text{Ad}U)(\mathcal D_n)$$
since $U^*U_1=D_U\in \mathcal D_n$. Moreover, we have:
\begin{equation}\begin{aligned}(\text{Ad}U_1)(e_2)&=(\text{Ad}U)\text{Ad}(\sum_{i,j=1}^n \bar u_{i,j}e_{j,j}\otimes e_{i,i})(\frac{1}{n}\sum_{k,l=1}^n e_{k,l})\\
&=(\text{Ad}U)(\sum_{i,k,l=1}^n \bar u_{i,k}u_{i,l}e_{k,l}\otimes e_{i,i})\\
&=(\text{Ad}U)(\text{Ad}U^*(e_3))\\
&=e_3
\end{aligned}\end{equation}

It follows that $AdU_1$ takes the basic construction $\mathbb C\subset\mathcal D_n\overset{e_2}\subset M_n(\mathbb C)$ onto the inclusion $\mathbb C\subset U\mathcal D_n U^*\overset{e_3}\subset U_1M_n(\mathbb C)U_1^*$. Thus this is also a basic construction, which shows that $\mathcal Q_1=U_1M_n(\mathbb C)U_1^*$. Moreover, it follows that each $\text{Ad}U_i$ takes the basic construction $\mathcal P_{i-1}\subset \mathcal P_i\subset \mathcal P_{i+1}$ onto $\mathcal Q_{i}\subset \mathcal Q_{i+1}\subset \mathcal Q_{i+2}$, which ends the proof.

\end{proof}

The first relative commutant $\mathcal D_n'\cap U\mathcal D_n U^*$ is equal to $\mathbb C$, since the commuting square condition implies $\mathcal D_n\cap U\mathcal D_n U^*=\mathbb C$. Thus the subfactor $N_H\subset M_H$ is irreducible. In the following proposition we describe the higher relative commutants of the subfactor $N_H\subset M_H$ as the commutants of some matrices $ P_i$, $i\geq 1$.

\begin{proposition} With the previous notations, let $P_i$ denote the projection $U_ie_{i+3}U_i^*\in \mathcal P_{i+1}$, $i\geq 1$. Then we have the following formula for the $(i+1)$-th relative commutant: $$\mathcal D_n'\cap \mathcal Q_i=P_i'\cap\mathcal D_n'\cap \mathcal P_i.$$

\end{proposition}
\begin{proof} We have: \begin{equation}\begin{aligned}
\mathcal D_n'\cap \mathcal Q_i&=\mathcal D_n'\cap\text{Ad}U_i(\mathcal P_{i-1})\\
&=\mathcal D_n'\cap\text{Ad}U_i(e_{i+3}'\cap\mathcal P_i)\\
&=\mathcal D_n'\cap P_i'\cap \text{Ad}U_i(\mathcal P_i)\\
&=\mathcal D_n'\cap P_i'\cap \mathcal P_i\\
\end{aligned}\end{equation}
We used the fact that $\mathcal P_{i-1}\subset \mathcal P_i\overset{e_{i+3}}\subset \mathcal P_{i+1}$ is a basic construction, and thus $e_{i+3}'\cap\mathcal P_i=\mathcal P_{i-1}$.
\end{proof}

\begin{remark} The $n^2\times n^2$ matrix $P_1=U_1e_4U_1^*$ can be written as $$P_1=\sum_{a,b,c,d=1}^n p_{a,b}^{c,d}e_{a,b}\otimes e_{c,d},\text{ where }p_{a,b}^{c,d}=\sum_{i=1}^n u_{a,i}\bar u_{b,i}\bar u_{c,i}u_{d,i}.$$
This matrix is used in the theory of Hadamard matrices and it is called the \emph{profile} of $H$. 
It is a result of Jones (\cite{Jo2}) that the matrices $P_{2i+1}$, $i\geq 1$, depend only on $P_1$. Indeed, one can check that $$P_{2i+1}=\sum_{k_1,l_1,...,k_i,l_i=1}^n p_{a,b}^{k_1,l_1}p_{k_1,l_1}^{k_2,l_2}...p_{k_i,l_i}^{c,d}e_{a,b}\otimes e_{k_1,l_1}\otimes e_{k_2,l_2}\otimes...\otimes e_{k_i,l_i}\otimes e_{c,d}.$$
Thus, all higher relative commutants of even orders are determined by $P_1$.
\end{remark}

Let $\Gamma_H$ denote the graph of vertices $\{1,2,...,n\}\times \{1,2,...,n\}$, in which the distinct vertices $(a,c)$ and $(b,d)$ are connected if and only if $p_{a,b}^{c,d}\not=0$. The second relative commutant can be easily described in terms of $\Gamma_H$. We recall this in the following Proposition, which is a reformulation of a result in \cite{Jo2} (see also \cite{JS}).
\begin{proposition}\label{second}The second relative commutant of the subfactor $N_H\subset M_H$ is abelian, its minimal projections are in bijection with the connected components of $\Gamma_H$, and their traces are proportional to the sizes of the connected components. 
\end{proposition}

\begin{proof} Let $\sum_{i,j=1}^n \lambda_i^j e_{i,i}\otimes e_{j,j}$, $\lambda_i^j\in\{0,1\}$,  be a projection in the second relative commutant $ P_1'\cap (\mathcal D_n\otimes \mathcal D_n)$. We have:
$$(\sum_{a,b,c,d=1}^n p_{a,b}^{c,d}e_{a,b}\otimes e_{c,d})(\sum_{i,j=1}^n \lambda_i^j e_{i,i}\otimes e_{j,j})=(\sum_{i,j=1}^n \lambda_i^j e_{i,i}\otimes e_{j,j})(\sum_{a,b,c,d=1}^n p_{a,b}^{c,d}e_{a,b}\otimes e_{c,d})$$
Equivalently:
$$\sum_{a,c,i,j=1}^n \lambda_i^j p_{q,i}^{c,j}e_{a,i}\otimes e_{c,j}=\sum_{b,d,i,j=1}^n \lambda_i^j p_{i,b}^{j,d} e_{i,b}\otimes e_{j,d}.$$
By relabeling and identifying the set of indices, it follows:
$$(\lambda_a^c-\lambda_i^j)p_{a,i}^{c,j}=0.$$
Thus, if the vertices $(a,c)$ and $(i,j)$ are connected then $\lambda_a^c=\lambda_i^j$. This ends the proof.
\end{proof}

\section{Matrices of Dita type}

In this section we investigate the standard invariant of subfactors associated to a particular class of Hadamard matrices, obtained by a construction of P.Dita (\cite{Di}), which is a generalization of an idea of U.Haagerup (\cite{Haagerup}). These matrices have a lot of symmetries, and we show that for such matrices the second relative commutant has some extra structure besides the Jones projection. 

Let $n$ be non-prime, $n=kl$ with $k,l\geq 2$. Let $A=(a_{i,j})\in M_{k}(\mathbb C)$ and $B_1,...,B_k\in M_{m}(\mathbb C)$ be complex Hadamard matrices. It is possible to construct an $n\times n$ Hadamard matrix from $A,B_1,...,B_k$ by using an idea of \cite{Di} (see also\cite{Haagerup},\cite{Petrescu}). This construction is a generalization of the tensor product of two Hadamard matrices:

\begin{equation}\label{Dita}H=\begin{pmatrix} a_{1,1}B_1 & a_{1,2} B_2 & ... & a_{1,k}B_k \\
a_{2,1}B_1 & a_{2,2} B_2 & ... & a_{2,k}B_k \\
.  & & & .\\
.  & & & .\\
.  & & & .\\
a_{k,1}B_1 & a_{k,2} B_2 & ... & a_{k,k}B_k \\
\end{pmatrix}\end{equation}

Let $(f_{i,j})_{1\leq i,j\leq k}$ be the matrix units of $M_{k}(\mathbb C)$. We identify $M_{n}(\mathbb C)$ with the tensor product $M_{m}(\mathbb C)\otimes M_{k}(\mathbb C)$, with the same conventions as before. Thus: $$H=\sum_{i,j=1}^k a_{i,j}B_j\otimes f_{i,j}$$

One can use construct multi-parametric families of non-equivalent Hadamard matrices, by replacing $B_1,...,B_k$ by $B_1D_1,...B_kD_k$, where $D_1,...,D_k$ are diagonal unitaries. Some of the families of Hadamard matrices of small orders considered in the next section arise from this construction.

Recall that the second relative commutant always contains the Jones projection $e_3=\sum e_{ii}\otimes e_{ii}$. In the next proposition we show that the second relative commutant of a Dita type subfactor contains another projection $f\geq e_3$, so it has dimension at least $3$. 

\begin{proposition}\label{Dita} Let $H=(a_{i,j}B_j)_{1\leq i,j\leq k}\in M_n(\mathbb C)$ be a Dita type matrix, where $A=(a_{i,j})_{1\leq i,j\leq k}\in M_k(\mathbb C)$ and $B_1,...,B_k\in M_m(\mathbb C)$ are complex Hadamard matrices, $n=mk$. Then the second relative commutant of the subfactor associated to $H$ contains the projection:
$$f=\sum_{1\leq i,j\leq n,\text{ } i\equiv j (\text{mod m)}} e_{i,i}\otimes e_{j,j}\in M_{n^2}(\mathbb C).$$
\end{proposition}

\begin{proof} For $1\leq i\leq n$, let $i_0=(i-1)\text{(mod m)}+1$ and $i_1=\frac{i-i_0}{m} +1$. We will use similar notations for $1\leq j \leq n$. Thus, the $(i,j)$ entry of $H$ is:
$$h_{i,j}=a_{i_1,j_1}b^{j_1}_{i_0,j_0}$$
where $b^t_{r,s}$ is the $(r,s)$ entry of $B_t$, for all $1\leq t\leq k$, $1\leq r,s\leq m$.

With these notations, the projection $f$ can be written as $$f=\sum_{i,j=1}^n \lambda_i^j e_{i,i}\otimes e_{j,j}$$
where $\lambda_i^j=1$ if $i_0=j_0$ and $\lambda_i^j=0$ for all other $i,j$. 

According to Proposition \ref{second}, showing that $f$ is in the second relative commutant is equivalent to showing that $p_{i,c}^{j,d}=0$ whenever $c_0\not= d_0$. Using the formula for the entries of $P_1$ and the fact that $i_0=j_0$ we obtain:
\begin{equation} \begin{aligned}p_{i,c}^{j,d}&=\sum_{x=1}^n u_{i,x}\bar u_{c,x}\bar u_{j,x} u_{d,x}\\
&=\frac{1}{n^2}\sum_{x=1}^n h_{i,x}\bar h_{c,x}\bar h_{j,x} h_{d,x}\\
&=\frac{1}{n^2}\sum_{x=1}^n a_{i_1,x_1}b^{x_1}_{i_0,x_0}\bar a_{c_1,x_1}\bar b^{x_1}_{c_0,x_0}\bar a_{j_1,x_1}\bar b^{x_1}_{j_0,x_0}a_{d_1,x_1}b^{x_1}_{d_0,x_0}\\
&=\frac{1}{n^2}\sum_{x=1}^n a_{i_1,x_1}\bar a_{c_1,x_1}\bar b^{x_1}_{c_0,x_0}\bar a_{j_1,x_1}a_{d_1,x_1}b^{x_1}_{d_0,x_0}\\
&=\frac{1}{n^2}\sum_{x_1=1}^k (a_{i_1,x_1}\bar a_{c_1,x_1}\bar a_{j_1,x_1}a_{d_1,x_1}(\sum_{x_0=1}^m\bar b^{x_1}_{c_0,x_0}b^{x_1}_{d_0,x_0}))\\
&=\frac{1}{n^2}\sum_{x_1=1}^k a_{i_1,x_1}\bar a_{c_1,x_1}\bar a_{j_1,x_1}a_{d_1,x_1}\delta_{c_0}^{d_0}\\
&=0\end{aligned}\end{equation}
whenever $c_0\not= d_0$.
\end{proof}

We show that in fact the subfactor $N_H\subset M_H$ associated to the Dita matrix $H$ has an intermediate subfactor $N_H\subset R_H\subset M_H$, and the projection $f$ is the Bisch projection (in the sense of \cite{Bi}) corresponding to $R_H$.

\begin{proposition}\label{adjacent} Let $H=\sum_{1\leq i,j\leq k} a_{i,j}B_j\otimes f_{i,j}\in M_n(\mathbb C)$ be a Dita type matrix, where $A=(a_{i,j})_{1\leq i,j\leq k}\in M_k(\mathbb C)$ and $B_1,...,B_k\in M_m(\mathbb C)$ are complex Hadamard matrices, $n=mk$. Then:

(a). The commuting square $\mathfrak{C}(H)$ can be decomposed into two adjacent symmetric commuting squares:
$$\begin{matrix}
\mathcal D_{m}\otimes D_{k} & \subset{} & M_m(\mathbb{C})\otimes M_k(\mathbb{C})\cr
\cr
\cup   &           & \cup \cr
\cr
\mathcal D_m\otimes I_k & \subset & U(M_m(\mathbb{C})\otimes \mathcal D_k)U^*\cr
\cr
\cup   &           & \cup \cr
\cr
\mathbb{C} & \subset{} & U\mathcal (\mathcal D_{m}\otimes D_{k})U^*
\end{matrix}$$

(b). The commuting square $\mathfrak{C}^t(H)$ can be decomposed into two adjacent symmetric commuting squares:
$$\begin{matrix}
U(\mathcal D_{m}\otimes D_{k})U^* & \subset{} & M_m(\mathbb{C})\otimes M_k(\mathbb{C})\cr
\cr
\cup   &           & \cup \cr
\cr
U(I_m\otimes \mathcal D_k)U^* & \subset & \mathcal D_m\otimes M_k(\mathbb{C})\cr
\cr
\cup   &           & \cup \cr
\cr
\mathbb{C} & \subset{} & \mathcal D_{m}\otimes D_{k}
\end{matrix}$$
\end{proposition}
\begin{proof}(a). We first show that $\mathcal D_m\otimes I_k \subset U(M_m(\mathbb{C})\otimes \mathcal D_k)U^*$. Equivalently, we check that $U^*(\mathcal D_m\otimes I_k)U\subset (M_m(\mathbb{C})\otimes \mathcal D_k)$. Indeed, for $D\in \mathcal D_m$ we have:
\begin{equation}\begin{aligned}
U^*(D\otimes I_k)U&=\frac{1}{n}(\sum_{1\leq i',j'\leq k} \bar a_{i',j'}B_{j'}^*\otimes f_{j',i'})(D\otimes I_k)(\sum_{1\leq i,j\leq k} a_{i,j}B_j\otimes f_{i,j})\\
&=\frac{1}{n}\sum_{1\leq i,j,j'\leq k}\bar a_{i,j'}a_{i,j}B_{j'}^*D B_j\otimes f_{j',j}\\
&=\frac{1}{n}\sum_{1\leq j,j'\leq k}(\sum_{i=1}^k\bar a_{i,j'}a_{i,j})B_{j'}^*D B_j\otimes f_{j',j}\\
&=\frac{1}{n}\sum_{1\leq j,j'\leq k}\delta_j^{j'}B_{j'}^*D B_j\otimes f_{j',j}\\
&=\frac{1}{n}\sum_{1\leq j\leq k}B_{j}^*D B_j\otimes f_{j,j}\in (M_m(\mathbb{C})\otimes \mathcal D_k)
\end{aligned}\end{equation}

The lower square of inclusions is clearly a commuting square, since $\mathfrak{C}(H)$ is a commuting square. We check that  
$$\begin{matrix}
\mathcal D_{m}\otimes D_{k} & \subset{} & M_m(\mathbb{C})\otimes M_k(\mathbb{C})\cr
\cr
\cup   &           & \cup \cr
\cr
\mathcal D_m\otimes I_k & \subset & U(M_m(\mathbb{C})\otimes \mathcal D_k)U^*\cr
\end{matrix}$$
is a commuting square. For $X\in M_m(\mathbb{C})$ and $D\in \mathcal D_k$ we have:
\begin{equation}\begin{aligned}
U(X\otimes D)U^*&=\frac{1}{n}(\sum_{1\leq i,j\leq k} a_{i,j}B_j\otimes f_{i,j})(X\otimes D)(\sum_{1\leq i',j'\leq k} \bar a_{i',j'}B_{j'}^*\otimes f_{j',i'})\\
&=\frac{1}{n}\sum_{1\leq i,i',j\leq k}\bar a_{i',j}a_{i,j}B_{j}X B_j^*\otimes D_{j,j}f_{i,i'}\\
\end{aligned}\end{equation}
Hence:
\begin{equation}\begin{aligned}
E_{\mathcal D_n}(U(X\otimes D)U^*)&=E_{\mathcal D_n}(\frac{1}{n}\sum_{1\leq i,i',j\leq k}\bar a_{i',j}a_{i,j}B_{j}X B_j^*\otimes D_{j,j}f_{i,i'})\\
&=\frac{1}{n}\sum_{1\leq i,i',j\leq k}E_{\mathcal D_m}(\bar a_{i',j}a_{i,j}B_{j}X B_j^*)\otimes D_{j,j}\delta_i^{i'}f_{i,i}\\
&=\frac{1}{n}\sum_{1\leq i,j\leq k}D_{j,j}E_{\mathcal D_m}(B_{j}X B_j^*)\otimes f_{i,i}\\
&=\frac{1}{n}\sum_{1\leq j\leq k}D_{j,j}E_{\mathcal D_m}(B_{j}X B_j^*)\otimes I_k\in \mathcal D_m\otimes \mathcal I_k
\end{aligned}\end{equation}
The lower commuting square is symmetric, since the product of the dimensions of its upper left and lower right corners equals the dimension of its upper right corner. This also implies that the upper commuting square is symmetric, since $\mathfrak{C}(H)$ is symmetric.

(b). The proof is similar to the proof of part (a).\end{proof}
\begin{corollary}The subfactors associated to Dita matrices have intermediate subfactors.
\end{corollary}
\begin{proof}By iterating the basic construction for the decomposition of $\mathfrak{C}^t(H)$ in commuting squares, we obtain the towers of algebras:
$$\begin{matrix}
U(\mathcal D_{m}\otimes D_{k})U^* & \subset{} & M_m(\mathbb{C})\otimes M_k(\mathbb{C})& \overset{e_3}\subset & \mathcal X_1 & \overset{e_4}\subset & \mathcal X_2 & \overset{e_5}\subset &...\cr
\cr
\cup   &           & \cup & & \cup & &\cup & \cr
\cr
U(I_m\otimes \mathcal D_k)U^*& \subset &  \mathcal D_m\otimes M_k(\mathbb{C})& \overset{e_3}\subset & \mathcal R_1 & \overset{e_4}\subset & \mathcal R_2 & \overset{e_5}\subset &...\cr
\cr
\cup   &           & \cup & & \cup & &\cup & \cr
\cr
\mathbb{C} & \subset{} & \mathcal D_{m}\otimes D_{k}& \overset{e_3}\subset & \mathcal Y_1 & \overset{e_4}\subset & \mathcal Y_2 & \overset{e_5}\subset &...
\end{matrix}$$
where $\mathcal R_i=<\mathcal R_{i-1},e_{i+2}>\subset \mathcal X_i$. Let $R_H$ be the weak closure of $\cup_i \mathcal R_i$. We have $N_H\subset R_H\subset M_H$ and $R_H$ is a $II_1$ factor since the subfactor $N_H\subset M_H$ is irreducible.

\end{proof}
\begin{remark} It is immediate to check that the projection $f\in M_n(\mathbb C)\otimes M_n(\mathbb C)$ from Proposition \ref{Dita} implements the conditional expectation from $M_n(\mathbb{C})\otimes I_n=M_n(\mathbb{C})$ onto $D_m\otimes M_k(\mathbb{C})$. It follows that $f$ is the Bisch projection for the intermediate subfactor $N_{H}\subset R_{H}\subset M_{H}$.
\end{remark}

\section{Matrices of small order}

In this section we compute the second relative commutants of the subfactors associated to Hadamard matrices of small dimensions. For some of the matrices considered we also specify the dimension of the third relative commutant. Most computations included were done with the help of computers, using GAP and Mathematica. 

It is well known in subfactor theory that the dimension of the second relative commutant $D'\cap Q_1$ is at most $n$, with equality if and only if $H$ is equivalent to a tensor product of Fourier matrices. In this case the subfactor $N_H\subset M_H$ is well understood, being a cross-product subfactor. For this reason, we exclude from our analysis tensor products of Fourier matrices.

Some of the matrices we present are parameterized and they yield continuous families of complex Hadamard matrices. In such cases, the strategy for computing the second relative commutant will be to determine which entries of the profile matrix $P_1$ depend on the parameters, and for what values of the parameters are these entries $0$. According to Proposition \ref{second}, the second relative commutant will not change as long as the $0$ entries of $P_1$ do not change. Thus, to compute the second relative commutant for any other value of the parameters, it is enough to compute it for some random value.

We will describe the second relative commutant by specifying its minimal projections. Each such projection $p$ corresponds to a subset $S\subset\{1,2,...,n^2 \}$: $p$ is the $n^2\times n^2$ diagonal matrix having $1$ on diagonal positions $i\in S$ and $0$ on all other positions. Since the Jones projection $e_3$ is always in the second relative commutant, one of the subsets of our partitions will always be $\{1,n+2,2n+3,...,kn+k+1,...,n^2\}$.

$$ $$

\textbf{Complex Hadamard matrices of dimension 4.} There exists, up to equivalence, only one family of complex Hadamard matrices of dimension 4:
$$F_{4}(a)=\begin{pmatrix} 1 & 1 & 1 & 1 \cr 1 & a & -1 & -a \cr 1 & 
    -1 & 1 & -1 \cr 1 & -a & -1 & a \cr  \end{pmatrix},\text{ }|a|=1$$

The entries of $P_1$ that depend on the parameter $a$ are $\frac{1}{8} + \frac{a^2}{8}$,
 $\frac{1}{8} - \frac{a^2}{8}$, $\frac{1}{8} + \frac{1}{8\,a^2}$, $\frac{1}{8} - \frac{1}{8\,a^2}$. Thus, the second relative commutant is the same for all values of $a$ that are not roots of these equations.

The roots $a=1, a=-1$ yield matrices that are tensor products of $2\times 2$ Fourier matrices. Thus the dimension of the second relative commutant is $4$, and its minimal projections are given by the partition $\{1,6,11,16\}$, $\{2,5,12,15\}$, $\{3,8,9,14\}$, $  \{4,7,10,13\}$.

The roots $a=i, a=-i$ yield the $4\times 4$ Fourier matrix, thus the minimal projections are $\{1,6,11,16\}$, $\{2,7,12,13\}$, $\{3,8,9,14\}$, $\{4,5,10,15\}$.

Any other values of $a$, $|a|=1$, yield relative commutants of dimension $3$: $\{1,6,11,16\}$, $\{2,4,5,7,10,12,13,15\}$, $\{3,8,9,14\} $. This is not surprising, since this matrix is of Dita type (see Proposition \ref{Dita}).

The dimension of the third relative commutant is $10$, and the dimension of the fourth relative commutant is $35$ unless $a$ is a primitive root of order $8$ of unity, in which case the dimension is $36$. Based on this evidence, we conjecture that the principal graph of the subfactor associated to $F_4(a)$ is $D_{2k}^{(1)}$ if $a$ is a primitive root of order $2^k$ of unity, and $D_{\infty}^{(1)}$ otherwise. 


$$ $$

\textbf{Complex Hadamard matrices of dimension 6.} The Fourier matrix $F_6$ is part of an affine 2-parameter family of Dita matrices:
$$F_{6}(a,b)=\begin{pmatrix} 1 & 1 & 1 & 1 & 1 & 1 \cr 1 & a\,e^{\frac{i }{3}\,\pi } & b\,e^{\frac{2\,i }{3}\,\pi } & -1 & \frac{a}
   {e^{\frac{2\,i }{3}\,\pi }} & \frac{b}{e^{\frac{i }{3}\,\pi }} \cr 1 & e^{\frac{2\,i }{3}\,\pi } & e^
   {\frac{-2\,i }{3}\,\pi } & 1 & e^{\frac{2\,i }{3}\,\pi } & e^{\frac{-2\,i }{3}\,\pi } \cr 1 & -a & b & 
    -1 & a & -b \cr 1 & e^{\frac{-2\,i }{3}\,\pi } & e^{\frac{2\,i }{3}\,\pi } & 1 & e^
   {\frac{-2\,i }{3}\,\pi } & e^{\frac{2\,i }{3}\,\pi } \cr 1 & \frac{a}{e^{\frac{i }{3}\,\pi }} & \frac{b}
   {e^{\frac{2\,i }{3}\,\pi }} & -1 & a\,e^{\frac{2\,i }{3}\,\pi } & b\,e^{\frac{i }{3}\,\pi } \cr \end{pmatrix}
$$

The entries of $P_1$ that depend on $a,b$ are: $2\,\left( 1 + a^{-2} + b^{-2} \right)$, $2 + \frac{2\,{\left( -1 \right) }^{\frac{2}{3}}}{a^2} - \frac{2\,{\left( -1 \right) }^{\frac{1}{3}}}{b^2}$, $2 - \frac{2\,{\left( -1 \right) }^{\frac{1}{3}}}{a^2} + \frac{2\,{\left( -1 \right) }^{\frac{2}{3}}}{b^2}$, $2\,\left( 1 + a^2 + b^2 \right)$, $2 + 2\,{\left( -1 \right) }^{\frac{2}{3}}\,a^2 - 2\,{\left( -1 \right) }^{\frac{1}{3}}\,b^2$, $2 - 2\,{\left( -1 \right) }^{\frac{1}{3}}\,a^2 + 2\,{\left( -1 \right) }^{\frac{2}{3}}\,b^2$.

Making one of these entries $0$ yields the following possibilities: $a=-\frac{1}{2} - \frac{i }{2}\,{\sqrt{3}}, b=-\frac{1}{2} + \frac{i }{2}\,{\sqrt{3}}$ or $a=-\frac{1}{2} + \frac{i }{2}\,{\sqrt{3}}, b=-\frac{1}{2} - \frac{i }{2}\,{\sqrt{3}}$ or $a=\frac{1}{2} - \frac{i }{2}\,{\sqrt{3}}, b=\frac{1}{2} + \frac{i }{2}\,{\sqrt{3}}$ or $a=\frac{1}{2} + \frac{i }{2}\,{\sqrt{3}}, b=\frac{1}{2} - \frac{i }{2}\,{\sqrt{3}}$ or $a=-\frac{1}{2} - \frac{i }{2}\,{\sqrt{3}}, b=\frac{1}{2} - \frac{i }{2}\,{\sqrt{3}}$ or $a=-\frac{1}{2} + \frac{i }{2}\,{\sqrt{3}}, b=\frac{1}{2} + \frac{i }{2}\,{\sqrt{3}}$ or $a=\frac{1}{2} - \frac{i }{2}\,{\sqrt{3}}, b=-\frac{1}{2} - \frac{i }{2}\,{\sqrt{3}}$ or $a=\frac{1}{2} + \frac{i }{2}\,{\sqrt{3}}, b=-\frac{1}{2} + \frac{i }{2}\,{\sqrt{3}}$ or $a=-1, b=-1$ or $a=1, b=1$ or $a=-1, b=1$ or $a=1, b=-1$.

In each of these cases the matrix $F_6(a,b)$ is a tensor product of Fourier matrices.

For all other pairs $(a,b)$ satisfying $|a|=|b|=1$, the second relative commutant has dimension $4$: $\{ 1,8,15,22,29,36\}$, $\{ 2,4,6,7,9,11,14,16,18,19,21,23,26,$\\$28,30,31,33,
   35\}$, $\{ 3,10,17,24,25,32\}$, $\{ 5,12,13,20,27,34\}$.

$ $

The following family of self-adjoint, non-affine, complex Hadamard matrices was obtained in \cite{BeN}, one of the motivations being the search for Hadamard matrices of small dimensions that might yield subfactors with no extra structure in their relative commutants, besides the Jones projections. 

$$BN_6(\theta)=\begin{pmatrix}
1 & 1 & 1 & 1 & 1 & 1 \cr
1 & -1 & \bar{x} &-y  &-\bar{x}  &y   \cr
1 & x & -1 & t & -t  & -x  \cr
1 &-\bar{y}  & \bar{t} & -1  & \bar{y} &-\bar{t}   \cr
1 &-x  & -\bar{t} & y & 1 &\bar{z}  \cr
1 & \bar{y} &-\bar{x}  &-t  & z & 1 \cr
\end{pmatrix}$$
where $\theta\in [-\pi,-arcos(\frac{-1+\sqrt{3}}{2})]\cup [arcos(\frac{-1+\sqrt{3}}{2}),\pi]$ and the variables $x,y,z,t$ are given by: 
$$y=exp(i\theta),\text{ } z=\frac{1+2y-y^2}{y(-1+2y+y^2)}$$
$$x=\frac{1+2y+y^2-\sqrt{2}\sqrt{1+2y+2y^3+y^4}}{1+2y-y^2}$$
$$t=\frac{1+2y+y^2-\sqrt{2}\sqrt{1+2y+2y^3+y^4}}{-1+2y+y^2}$$
The entries of $BN_6$ do not depend linearly on the parameters, thus this is not a Dita-type family. The corresponding subfactors have the second relative commutant generated by the Jones projection. We conjecture that $BN_6(\theta)$ give supertransitive subfactors, i.e. all the relative commutants of higher orders are generated by the Jones projections.

There are other interesting complex Hadamard matrices of order $6$, such as the one found by Tao in connection to Fuglede's conjecture (\cite{T}), or the Haagerup matrix (\cite{Haagerup},{TZ}). We computed the second relative commutant for these matrices, and it is generated by the Jones projection.

$ $

\textbf{Complex Hadamard matrices of dimension 7.} The following one-parameter family was found in \cite{Petrescu}, providing a counterexample to a conjecture of Popa regarding the finiteness of the number of complex Hadamard matrices of prime dimension.

$${P_{7}}(a)=\begin{pmatrix} 1 & 1 & 1 & 1 & 1 & 1 & 1 \cr 1 & a\,
   e^{\frac{i }{3}\,\pi } & \frac{a}
   {e^{\frac{2\,i }{3}\,\pi }} & e^{\frac{-i }{3}\,\pi } & 
    -1 & -1 & e^{\frac{i }{3}\,\pi } \cr 1 & \frac{a}
   {e^{\frac{2\,i }{3}\,\pi }} & a\,e^{\frac{i }{3}\,\pi } & 
    -1 & e^{\frac{-i }{3}\,\pi } & -1 & e^
   {\frac{i }{3}\,\pi } \cr 1 & e^{\frac{-i }{3}\,\pi } & 
    -1 & \frac{e^{\frac{i }{3}\,\pi }}{a} & \frac{1}
   {a\,e^{\frac{2\,i }{3}\,\pi }} & e^{\frac{i }{3}\,\pi } & 
    -1 \cr 1 & -1 & e^{\frac{-i }{3}\,\pi } & \frac{1}
   {a\,e^{\frac{2\,i }{3}\,\pi }} & \frac{e^
     {\frac{i }{3}\,\pi }}{a} & e^{\frac{i }{3}\,\pi } & 
    -1 \cr 1 & -1 & -1 & e^{\frac{i }{3}\,\pi } & e^
   {\frac{i }{3}\,\pi } & e^{\frac{-2\,i }{3}\,\pi } & e^
   {\frac{-i }{3}\,\pi } \cr 1 & e^{\frac{i }{3}\,\pi } & e^
   {\frac{i }{3}\,\pi } & -1 & -1 & e^
   {\frac{-i }{3}\,\pi } & e^{\frac{-2\,i }{3}\,\pi } \cr  \end{pmatrix}$$

The second relative commutant of the associated subfactors is generated by the Jones projection, for all $|a|=1$. For $a=1$ we also computed the third relative commutant, and it is just the Temperley-Lieb algebra $TL_2$. We conjecture that $P_7(a)$ yield subfactors with no extra structure in their higher order relative commutants, besides the Jones projections.

$ $

\textbf{Complex Hadamard matrices of dimension 8.} The following 5-parameter family of Hadamard matrices contains the Fourier matrix and is of Dita type:

$$F_{8}(a,b,c,d,z)=\begin{pmatrix} 1 & 1 & 1 & 1 & 1 & 1 & 1 & 1 \cr 1 & a\,
   e^{\frac{i }{4}\,\pi } & i \,b & c\,
   e^{\frac{3\,i }{4}\,\pi } & -1 & \frac{a}
   {e^{\frac{3\,i }{4}\,\pi }} & -i \,b & \frac{c}
   {e^{\frac{i }{4}\,\pi }} \cr 1 & i \,d & -1 & -i \,
   d & 1 & i \,d & -1 & -i \,d \cr 1 & e^
    {\frac{3\,i }{4}\,\pi }\,z & -i \,b & \frac{c\,
     e^{\frac{i }{4}\,\pi }\,z}{a} & -1 & \frac{z}
   {e^{\frac{i }{4}\,\pi }} & i \,b & \frac{c\,z}
   {a\,e^{\frac{3\,i }{4}\,\pi }} \cr 1 & -1 & 1 & -1 & 1 & 
    -1 & 1 & -1 \cr 1 & \frac{a}{e^{\frac{3\,i }{4}\,\pi }} & 
   i \,b & \frac{c}{e^{\frac{i }{4}\,\pi }} & -1 & a\,
   e^{\frac{i }{4}\,\pi } & -i \,b & c\,
   e^{\frac{3\,i }{4}\,\pi } \cr 1 & -i \,d & -1 & i \,
   d & 1 & -i \,d & -1 & i \,d \cr 1 & \frac{z}
   {e^{\frac{i }{4}\,\pi }} & -i \,b & \frac{c\,z}
   {a\,e^{\frac{3\,i }{4}\,\pi }} & -1 & e^
    {\frac{3\,i }{4}\,\pi }\,z & i \,b & \frac{c\,
     e^{\frac{i }{4}\,\pi }\,z}{a} \cr  \end{pmatrix}$$

The list of possible values of $a,b,c,d,z$ that yield $0$ entries for $P_1$ is very long and we do not include it here. Outside these values, the second relative commutant has dimension $4$ and it is given by $\{1,10,19,28,37,46,55,64\}$, $\{2,4,6,8,9,11,13,15,18,20,22,24,25,27,29,31,34,36,
    38,40,41,43,45,47,50,$\\$52,54,56,57,59,61,63\}$, $\{3,7,12,16,17,21,26,30,35,39,
    44,48,49,53,58,62\}$, $\{5,14,23,32,33,42,51,60\}$. 

$ $


%

$ $

We analysed several other complex Hadamard matrices besides those included in this paper, such as those found by \cite{MRS},\cite{Sz}. We tried to cover most known examples of complex Hadamard matrices of dimensions $2,3,...,11$. We draw some conclusions:
\begin{enumerate}
\item Matrices of Dita type yield subfactors with intermediate subfactors, and thus the second relative commutant has some extra structure besides the Jones projection. We note that parametric families of Dita type exist for every $n$ non-prime, and they contain the Fourier matrix $F_n$.

\item All non-Dita, non-Fourier matrices we tested have the second relative commutant generated by the Jones projection. The third relative commutant is also generated by the first two Jones projections for all cases we could compute. It remains an open problem whether there exist such complex Hadamard matrices that admit symmetries of higher order.

\end{enumerate}

\bibliographystyle{amsalpha}

\end{document}